\documentclass[12pt]{amsart}
\usepackage{amssymb, amstext, amscd, amsmath}
\makeatletter
\def\@cite#1#2{{\m@th\upshape\bfseries%
[{#1\if@tempswa{\m@th\upshape\mdseries, #2}\fi}]}}
\makeatother
\theoremstyle{plain}
\newtheorem{thm}{Theorem}[section]
\newtheorem{cor}[thm]{Corollary}
\newtheorem{prop}[thm]{Proposition}
\newtheorem{lem}[thm]{Lemma}
\theoremstyle{definition}
\newtheorem{rem}[thm]{Remark}
\newtheorem{defn}[thm]{Definition}

\newtheorem{eg}[thm]{Example}
\newcommand{\Prf}{\noindent\textbf{Proof.\ }}
\newcommand{\bx}{\hfill$\blacksquare$\medbreak}
\newcommand{\upbx}{\vspace{-2.5\baselineskip}\newline\hbox{}%
\hfill$\blacksquare$\newline\medbreak}

\newcommand{\ca}{\mathrm{C}^*}

\newcommand{\bbC}{{\mathbb{C}}}

\newcommand{\bbD}{{\mathbb{D}}}
\newcommand{\bbK}{{\mathbb{K}}}
\newcommand{\bbN}{{\mathbb{N}}}

\newcommand{\bbR}{{\mathbb{R}}}
\newcommand{\bbT}{{\mathbb{T}}}
\newcommand{\bbZ}{{\mathbb{Z}}}

  \newcommand{\A}{{\mathcal{A}}}
  \newcommand{\B}{{\mathcal{B}}}

\renewcommand{\H}{{\mathcal{H}}}
  
  \newcommand{\J}{{\mathcal{J}}}
  \newcommand{\K}{{\mathcal{K}}}

  \newcommand{\M}{{\mathcal{M}}}

\renewcommand{\P}{{\mathcal{P}}}

\renewcommand{\S}{{\mathcal{S}}}
  
  \newcommand{\U}{{\mathcal{U}}}

  \newcommand{\X}{{\mathcal{X}}}

\newcommand{\eps}{\varepsilon}
\newcommand{\fA}{{\mathfrak{A}}}

\newcommand{\fT}{{\mathfrak{T}}}

\newcommand{\AND}{\text{ and }}
\newcommand{\FOR}{\text{ for }}

\newcommand{\qand}{\quad\text{and}\quad}

\newcommand{\qfor}{\quad\text{for}\quad}
\newcommand{\qforal}{\quad\text{for all}\quad}

\newcommand{\alg}{\operatorname{alg}}

\newcommand{\id}{{\operatorname{id}}}

\newcommand{\spn}{\operatorname{span}}

\newcommand{\rep}{\operatorname{rep}}

\newcommand{\bz}{\mathbf{z}}

\begin{document}

%%%%%%%%%%%%%%%%%%%%%%%%%%%%%%%%%%%%%%%%%%%%%
\title[Conjugacy Algebras]{Isomorphisms between Topological Conjugacy Algebras.}
%\thanks{}
%

%
\author[K.R.Davidson]{Kenneth R. Davidson}
\thanks{First author partially supported by an NSERC grant.
Elias Katsoulis was partially supported by
a summer grant from ECU}
\address{Pure Math.\ Dept.\\U. Waterloo\\Waterloo, ON\;
N2L--3G1\\CANADA}
\email{krdavids@uwaterloo.ca}

\author[E.G. Katsoulis]{Elias~G.~Katsoulis}
%\thanks{Katsoulis' research partially supported by a grant from
%ECU}
\address{Department of Mathematics\\East Carolina University\\
Greenville, NC 27858\\USA}
\email{KatsoulisE@mail.ecu.edu}
\begin{abstract}
A family of algebras, which we call topological conjugacy
algebras, is associated with each proper continuous map
on a locally compact Hausdorff space.
Assume that $\eta_i:\X_i\to \X_i$ is a
continuous proper map on a locally compact Hausdorff space $\X_i$,
for $i = 1,2$. We show that the dynamical systems $(\X_1, \eta_1)$
and $(\X_2, \eta_2)$ are conjugate
if and only if some topological conjugacy algebra of $(\X_1, \eta_1)$
is isomorphic as an algebra to
some topological conjugacy algebra of $(\X_2, \eta_2)$.
This implies as a corollary the complete classification of the
semicrossed products $C_0(\X) \times_{\eta} \bbZ^{+}$,
which was previously considered by Arveson and Josephson \cite{ArvJ},
Peters \cite{Pet}, Hadwin and Hoover \cite{HadH} and Power \cite{Pow}.
We also obtain a complete classification of all semicrossed products of the form
$A(\bbD) \times_{\eta}\bbZ^{+}$, where $A(\bbD)$ denotes the
disc algebra and $\eta: \bbD \rightarrow \bbD$ a continuous map
which is analytic on the interior.
In this case, a surprising dichotomy appears in the classification scheme,
which depends on the fixed point set of $\eta$.
We also classify more general semicrossed products of uniform algebras.
\end{abstract}

\date{}
\maketitle
%%%%%%%%%%%%%%%%%%%%%%%%%%%%%%%%%%%%%%%%%%%%%%%
\section{introduction}

The main objective of this paper is the complete classification of various classes of semicrossed
products that have appeared in the literature through the years.
One such class consists of Peters' semicrossed products
$C_0(\X) \times_{\eta}\bbZ^+$, where $(\X, \eta)$ is a dynamical system consisting of a proper continuous map $\eta$ acting on a
locally compact Hausdorff space. Under the assumption that
the topological spaces are compact and the maps are aperiodic, Peters \cite{Pet}
showed that two such semicrossed products are isomorphic
as algebras if and only if the corresponding dynamical systems are conjugate,
thus extending an earlier classification scheme of Arveson \cite{Arv} and
Arveson and Josephson \cite{ArvJ}. In \cite{HadH}
Hadwin and Hoover considered more general dynamical systems
of the form $(\X, \eta)$, where $\X$ is a compact Hausdorff space and $\eta$ is a continuous map so that
the set  $\{ x \in \X : \eta(x) \neq x, \, \eta^{(2)}(x)=\eta(x)\}$ has empty interior. For such
dynamical systems, it was shown \cite[Theorem 3.1]{HadH} that the Arveson-Josephson-Peters classification scheme holds.
In particular, \cite[Theorem 3.1]{HadH} classified all semicrossed products of the form
$C_0(\X) \times_{\eta}\bbZ^+$,
where $\X$ is compact Hausdorff and $\eta$ a homeomorphism. The analogous result for locally compact Hausdorff spaces
was proven by Power in \cite{Pow}. (Power's technique is applicable to a more general setting, which is limited
however by assumptions similar to those of Hadwin and Hoover regarding the fixed point set of $\eta$.) In spite of the continuing interest on semicrossed products
and their
variants \cite{AlaP, Bu, BP, DKM, HPW, Lam-1, Lam, MM, OS, Pow, Pow2, Solel-1}, the problem of classifying semicrossed products of the form
$C_0(\X) \times_{\eta}\bbZ^+$ remained open in the generality introduced by
Peters in \cite{Pet}. This problem is now being resolved in this paper by showing that the
Arveson-Josephson-Peters classification scheme holds with no restrictions on either $\X$
or $\eta$.

In this paper, we view semicrossed products of the form
$C_0(\X) \times_{\eta}\bbZ^+$ as a special case of a class
of Banach algebras that we call \textit{topological
conjugacy algebras} (see Definitions \ref{conjalg1} and \ref{conjalg2}).
The reader familiar with the work of Hadwin and Hoover \cite{HadH}
will recognize that in the compact case, our
topological conjugacy algebras are precisely the algebras that satisfy the
requirements of \cite[Proposition 2.7]{HadH}. Therefore in that case, our topological conjugacy
algebras form a proper subclass of the \textit{conjugacy algebras} of Hadwin and Hoover.
However, all known examples of the Hadwin-Hoover algebras which are
Banach algebras obey our definitions; and it seems that
little generality is lost by focusing on our class.
Most importantly, our context is applicable to locally compact spaces as well.
We develop a theory of conjugacy algebras for pairs $(\X, \eta)$, where $\eta$
is a proper continuous map on a locally compact Hausdorff space.
(Such a theory was promised in \cite{HadH} but was never delivered.)
One of the main results of the paper, Theorem \ref{main}, shows that
two dynamical systems $(\X_1, \eta_1)$ and $(\X_2, \eta_2)$ are conjugate
if and only if some topological conjugacy algebra of $(\X_1, \eta_1)$ is isomorphic
as an algebra to some topological conjugacy algebra of $(\X_2, \eta_2)$.
This result provides the desired classification of
Peters' semicrossed products discussed earlier, and also gives an affirmative answer
to the conjecture raised by Hadwin and Hoover in \cite[Remark 3.1.(2)]{HadH}.

A class of algebras related to Peter's semicrossed products \break $C_0(\X) \times_{\eta}\bbZ^+$
are the semicrossed products of the disc algebra,
i.e., semicrossed products of the form $A(\bbD) \times_{\eta}\bbZ^{+}$, where $A(\bbD)$ denotes the
disc algebra and $\eta: \bbD \rightarrow \bbD$ a continuous map which is analytic in the interior.
These algebras have been studied in \cite{Bu, BP, Hoov, HPW, Pow2} and for
the last fifteen years they
have been the source of many interesting theorems and counterexamples in
operator theory. Most notable is the dilation result of Buske and Peters \cite[Theorem II.4]{BP}
 (see also \cite{Pow2}). This result is a non-commutative analogue of Ando's Theorem and
 shows that for a conformal mapping $\eta$,
 the semicrossed product $A(\bbD) \times_{\eta}\bbZ^{+}$ plays the same role
 for a pair of $\eta$-commuting contractions
that the disc algebra $A(\bbD)$ plays for a single contraction. Since
$A(\bbD) \times_{\eta}\bbZ^{+}$ is a universal object,
it provides a non-commutative functional calculus (and therefore von Neumann type
inequalities) for pairs of $\eta$-commuting contractions.

In spite of the importance of the semicrossed products of the disc algebra,
 very little is available regarding their
isomorphic classes. (This lack of knowledge is apparent in
\cite[Remark I.V]{BP}, where the authors identify only three different isomorphic classes.)
By elaborating on the techniques of Theorem \ref{main}, we now obtain a complete classification
of these algebras. It turns out that there is a dichotomy in their
classification scheme, depending on the structure of the fixed point set of $\eta$.
If $\eta$ is not elliptic, then the (algebraic) isomorphism class of
$A(\bbD) \times_{\eta}\bbZ^{+}$ is determined by the conjugacy class of $\eta$,
with respect to conformal maps of the open disc (Corollary \ref{maindisc1}).
On the other hand, if $\eta$ is elliptic, then the isomorphism class of
$A(\bbD) \times_{\eta}\bbZ^{+}$ is determined
by the conjugacy classes of both $\eta$ and $\eta^{-1}$ (Theorem \ref{maindisc2}).
Surprisingly, this dichotomy disappears if one considers semicrossed products
of the form $A(\bbK) \times_{\eta}\bbZ^{+}$, where $\bbK$ is the closure of a
bounded region which is not simply connected and
and its boundary consists of finitely many Jordan curves (e.g. an annulus).
In that case, the isomorphism class of $A(\bbK) \times_{\eta}\bbZ^{+}$
is determined by the conjugacy class of $\eta$, with respect to conformal maps of $\bbK$
(Theorem \ref{generalregion}).

A final remark regarding the proofs. Even though
our initial approach is based on earlier ideas of others regarding
the character space of a conjugacy algebra, the ``hard'' part of the proofs
depend on a new idea, the use of two dimensional nest representations.
Such representations were first used in the study of classification
problems by the second author and Kribs \cite{KK} (see also \cite{Solel}).

\section{Definitions and examples}

Let $\X$ be a compact Hausdorff space and let $C(\X)$ denote the
continuous functions on $\X$. If $\eta:\X \rightarrow \X$ is
a continuous function on $\X$, then the \textit{skew polynomial algebra}
$P(\X, \eta)$ consists of all polynomials of the form
$\sum_n f_nU^n$, $f_n\in C(\X)$, where the multiplication of the
``coefficients'' $f\in C(\X)$ with the ``variable'' $U$ obeys the rule
\[
Uf = (f\circ \eta)  U
\]

Now consider
$\A$ to be a Banach algebra which satisfies:
\begin{itemize}
\item[(1)] the skew polynomial algebra $P(\X, \eta)$ is contained
as a dense subalgebra of $\A$, and the
constant function $1$ is the identity for $\A$.
\item[(2)] $C(\X) \subseteq P(\X, \eta) \subseteq \A$ is closed and there exists an algebra
homomorphism $E_0:\A \rightarrow C(\X)$ so that:
\begin{itemize}
\item[(i)] $E_0(f) = f$, for all $f \in C(\X)$, and
\item[(ii)] $\ker E_0 =  \A U$.
\end{itemize}
\item[(3)] $U$ is not a right divisor of $0$.
\end{itemize}

Since $C(\X)$ is closed in $\A$, a classical automatic continuity
result implies that $E_0$ is continuous and therefore $\A U =\ker E_0$ is closed as well.
Thus, the Inverse Mapping Theorem implies that the injection
$S:\A\to\A U$ by $Sa = aU$
has a bounded left inverse, which we denote by $T$.

With each $a \in \A$ we now associate a formal power series
\begin{equation} \label{powerseries}
a \sim \sum_n  E_n (a)U^n \in
 P^{\infty}(X, \eta).
 \end{equation}
  The map $E_0$ has already been defined above.
 Since
 \[
 a - E_0(a) \in \ker E_0 = \A U,
 \]
  there is a unique
 $b \in A$ so that $a = E_0(a) + b U$. We can now inductively define the maps
 $E_n, n=1, 2, \dots$. An easy computation shows that
$E_n= E_0(T(I-E_0))^n$ for $n \ge 1$, and in particular
the coefficient maps $E_n$ are bounded.
Then we can easily verify that the map
 \begin{equation} \label{Delta}
 \Delta : \A \longrightarrow P^{\infty}(X, \eta)
 \quad\text{given by} \quad \Delta(a) = \sum_{n\ge0} E_n (a) U^n
 \end{equation}
 is an algebra homomorphism.

\begin{defn} \label{conjalg1}
 Let $\X$ be a compact Hausdorff space, $\eta:\X \rightarrow \X$ be a
 continuous function and $\A$ be a Banach algebra that satisfies
 conditions (1), (2) and (3) above. Then $\A$ is said to
 be a \textit{topological conjugacy algebra} for $(\X, \eta)$ if
 $\limsup_n (\|E_n\|\|U^n\|)^{1/n} \leq 1$.
\end{defn}

In order to define topological conjugacy algebras for dynamical systems on a non-compact,
locally compact Hausdorff space $\X$,
we use the one-point compactification $\hat{X} = \X \cup \{\omega\}$ of such a space. In that case,
we will identify $C_0(\X)$ (the continuous functions on $\X$ ``vanishing at infinity''
with the continuous function on $\hat{\X}$ vanishing at $\omega$.
Any continuous proper map on $\X$ extends to a continuous map on
$\hat{\X}$ which has $\omega$ as a fixed point.

\begin{defn} \label{conjalg2}
Let $\X$ be a non-compact, locally compact Hausdorff space,
and let $\eta$ be a proper continuous map on $\X$.
Let $\hat{\A}$ be a topological conjugacy algebra for $(\hat{\X}, \hat{\eta})$.
Then, the norm closed algebra $\A$, which is
generated by the polynomials with coefficients in $C_0(\X)$ is said to be
a \textit{topological conjugacy algebra} for $(\X , \eta)$. Furthermore, the algebra $\hat{\A}$
is called the \textit{canonical unitization} of $\A$.
\end{defn}

For the rest of the paper we adopt the following conventions regarding our notation. Whenever we
say that a topological conjugacy algebra $\A$ for a dynamical system
$(\X, \eta)$ has a canonical unitization
$\hat{\A}$, we automatically imply that $\X$ is non-compact,
the one point compactification of $\X$ is denoted as $\hat{\X}$ and that
the extension of $\eta$ to $\hat{\X}$ is denoted as $\hat{\eta}$.

 Note that the formal power series expansion of (\ref{powerseries}) is applicable to conjugacy algebras defined on non-compact
 spaces as well. Indeed, let $\A$ be a conjugacy algebra for a dynamical system $(\X, \eta)$
with canonical unitization $\hat{\A}$. Since $\hat{\X}$ is compact, we have a power series expansion
 \[
 a \sim \sum_n  \hat{E}_n (a)U^n \in
 P^{\infty}(\hat{X}, \hat{\eta})
 \]
for any $a \in \hat{\A}$. Since $\hat{E}_n (P(\X , \eta))\subseteq C_0(\X)$, the continuity of
$\hat{E}_n$ implies that $\hat{E}_n (\A)\subseteq C_0(\X)\subseteq \A$, and the existence of the power
series expansion follows. The restriction of the coefficient maps $\hat{E}_n$ on $\A$ will
be denoted as $E_n$.

\begin{prop}  \label{Cohen}
Let $\X$ be a locally compact Hausdorff space and $\eta$ a proper continuous map on $\X$.
Let $\A$ be a topological conjugacy algebra for $(\X, \eta)$, $\B$ an algebra
and $\rho: \A \rightarrow \B$ an algebra homomorphism. If $C_0(\X)U \subseteq
\ker \rho$, then $\rho(a) = \rho(E_0(a))$, for all $a \in \A$.
\end{prop}

\Prf If $a \in \A$, then
\[
a = E_0(a)+ E_1(a)U + bU^2
\]
for some $b \in \A$. By our assumptions, $\rho$ annihilates the second summand in the above sum. The conclusion
will follow if we show that $\rho$ annihilates $bU^2$ as well.

Let $(e_{\alpha})$ be the contractive approximate unit of $C_0(\X)$
of all positive functions of norm less than $1$ with compact support. (The net $(e_{\alpha})$
is ordered with the usual pointwise order, i.e.,
$e_{\alpha}\leq e_{\beta}$ if and only if $e_{\alpha}(x)\leq e_{\beta}(x)$,
for all $x \in \X$.) Note that $(e_i)$ is also an approximate
unit for the Banach algebra $\A U^2$ and so $\overline{\spn}(\A U^2 C_0(\X))=\A U^2$.
Hence the right multiplication on $\A U^2$ by elements of
$C_0(\X)$ defines an anti-representation of $C_0(\X)$ on the Banach space $\A U^2$ which satisfies
the hypothesis of Cohen's Factorization Theorem \cite[Theorem 5.2.2]{Palm}. The conclusion is that
\[
\overline{\spn}(\A U^2 C_0(\X))= \spn(\A U^2 C_0(\X)) =\A U^2
\]
and so there exist $c_i \in \A$ and $g_i\in C_0(\X)$ so that
$bU^2  = \sum_{i=1}^n c_i U^2 g_i$.  But then
\[
\rho(bU^2) = \sum_{i=1}^n \rho(c_iU) \rho(Ug_i) =
  \sum_{i=1}^n \rho(c_iU) \rho((g_i \circ \eta)U) = 0 .
\]
and the conclusion follows.
\bx

We now present some examples of topological conjugacy algebras, including the semicrossed
products discussed in the introduction.

\begin{eg}   \label{rightone}
Let $\X$ be a locally compact space and $\eta$ a proper continuous map on $\X$.
Let $\H_x = l^2(\bbN)$ be the space of square summable sequences $\xi = (\xi_n)_{n=0}^{\infty}$
and for each $x \in \X$ define
\[
\pi_x(f)\xi = (f(x)\xi_0, (f\circ \eta)(x)\xi_1, (f\circ\eta^{(2)})(x)\xi_2, \dots),
\qfor f \in C_0(\X),
\]
and
\[
V_x\xi = (\xi_1, \xi_2, \xi_3, \dots).
\]
The norm closed operator
algebra $\A_{\X, \eta}$ acting on $\oplus_{x \in \X}\H_x$ and generated by the
operators,
\begin{equation*}
\oplus_{x \in \X} \pi_x(f),\quad
 \oplus_{x \in \X} \pi_x(g)V_x, \quad f,g \in C_0(\X),
\end{equation*}
is easily seen to be a topological conjugacy algebra.
\end{eg}

\begin{eg} \label{exsemicrossed}
(\textit{The semicrossed product $C_0(\X) \times_{\eta}\bbZ^+$} of \cite{Pet}) \\
Let $\X$, $\H_x$, and $\pi_x(f)$ be as above, and let
$U_x$ be the forward shift
\[
U_x\xi = (0, \xi_0, \xi_1, \xi_2, \dots).
\]
We define $C_0(\X) \times_{\eta}\bbZ^+$ to be the norm closed operator
algebra acting on $\oplus_{x \in \X}\H_x$ and generated by the
operators,
\begin{equation} \label{generators}
\oplus_{x \in \X} \pi_x(f), \quad
 \oplus_{x \in \X} U_x\pi_x(g), \quad f,g \in C_0(\X).
\end{equation}
This is not a conjugacy algebra in our sense because
the relation satisfied by these representations is
\[ \pi_x (f)U_x = U_x \pi_x((f\circ \eta)) \]
rather than the other way around.  But there is a natural connection.

Let $(C_0(\X) \times_{\eta}\bbZ^+)_{op}$ denote the opposite algebra of
$C_0(\X) \times_{\eta}\bbZ^+$, i.e.,
the algebra resulting from $C_0(\X) \times_{\eta}\bbZ^+$ by retaining
the same addition and introducing a new
multiplication $\circledcirc $ defined by $a \circledcirc b\equiv ba$
for $a, b \in C_0(\X) \times_{\eta}\bbZ^+$.
It is readily verified that $(C_0(\X) \times_{\eta}\bbZ^+)_{op}$ satisfies Definition
\ref{conjalg2} and is therefore a topological conjugacy algebra. In general,
Peters' semicrossed product $C_0(\X) \times_{\eta}\bbZ^+$ might not be isomorphic
to $(C_0(\X) \times_{\eta}\bbZ^+)_{op}$. Nevertheless, the semicrossed
products $C_0(\X) \times_{\eta_1}\bbZ^+$ and $C_0(\X) \times_{\eta_2}\bbZ^+$ are isomorphic
as algebras if and only if their opposite algebras are and so our results regarding isomorphisms between
topological conjugacy algebras have the analogous implications
for isomorphisms between Peters' semicrossed products.

In the case where $\eta$ is a surjection, Peters defines the semicrossed product
$C_0(\X) \times_{\eta}\bbZ^+$ as the universal operator algebra that makes any isometric covariant representations
of $(C_0(\X) , \eta)$ contractive on the skew polynomial algebra $P(\X, \eta)$ equiped with the $l^1$-norm \cite[Definition II.2]{Pet}. In
\cite[Proposition II.7]{Pet} it is shown that this universal object is isomorphic to the operator algebra
we have just defined. In the case where $\eta$ is not a surjection,
the semicrossed product $C_0(\X) \times_{\eta}\bbZ^+$ satisfies a similar universal property but to show this one
needs to use the technical language of
$\ca$-correspondences \cite{Kat}. We plan to pursue this elsewhere.

Finally, we remark that there are definitive relations between \break $C_0(\X) \times_{\eta}\bbZ^+$,
$(C_0(\X) \times_{\eta}\bbZ^+)_{op}$ and the operator algebras in Example \ref{rightone} but they
need not concern us here.
\end{eg}

\begin{eg} By letting the operators in (\ref{generators}) act on other
spaces of summable
sequences, we obtain additional examples of topological conjugacy algebras.
\end{eg}

The paper of Hadwin and Hoover
\cite{HadH} contains additional examples of conjugacy
algebras and the interested reader is directed there.

%%%%%%%%%%%%%%%%%%%%%%%%%%%%%%%%%%%%%%%%%%%%%%%%
\section{The character space of a conjugacy algebra}

Let $\X$ be a locally compact Hausdorff space, let $\eta$ a continuous map on $\X$ and $\A$ a
conjugacy algebra for $(\X, \eta)$. If $\rho$ is a character for $\A$, then its
action on $C_0(\X)\subseteq \A$ is a point evaluation on some point $x \in \X$.
This forces a partition of
the character space $\M_{\A}$ into disjoint sets $\M_{\A, x}$, $x \in \X$,
each one consisting of the
characters that coincide with the point evaluation on $x$ when restricted to
$C_0(\X)$.

If $x \in \X$ is
not a fixed point for $\eta$, then a simple argument with the skew relation
\[
g U h = (h\circ\eta) gU \qfor g, h \in C_0(\X),
\]
shows that $C_0(\X)U \subseteq
\ker \theta$ for any $\theta \in \M_{\A, x}$.
By Proposition \ref{Cohen}, $\theta(a) = \theta(E_0(a))= E_0(a)(x)$,
for all $a \in \A$ and so  $\M_{\A, x}$ contains a single element,
denoted as $\theta_{x,0}$.

Assume that $x$ is a fixed point for $\eta$ and $\X$ is compact.
In that case, Hadwin and Hoover \cite{HadH}
have shown that $\M_{\A, x}$ is homeomorphic to the spectrum
$\sigma (U)$ of $U$, which happens to be a \textit{closed} disc $\bbD_r$
centered at the origin with positive radius $r$.
We repeat their argument for completeness.

Indeed, let $r = \lim_{n\to\infty} \|U^n\|^{1/n}$. Recall that the map
$Sa = aU$ from $\A$ to $\A U$
has a bounded left inverse $T$; and is therefore bounded below by $\|T\|^{-1}$.
Hence, $\|U\|\geq \|T\|^{-1} \|1\|$ and so $\|U^n\|\geq \|T\|^{-n}\|1\|$.
Thus $r\geq \|T\|^{-1}$ is positive.

Now consider the power series $\sum_n E_n(a)(x)z^n$ and note that by Definition \ref{conjalg1},
its radius of convergence is at least $r$. Hence, for any $z \in \bbC$ with $|z|<r$, the mapping
\[
\A \ni a \longrightarrow  \sum_n E_n(a)(x)z^n
\]
is a well defined multiplicative functional on $\A$,
which we denote by $\theta_{x, z}$.
Since $\theta_{x,z}(U)=z$, for any $z$ with $|z|<r$, we have $\sigma(U)=\bbD_r$.
Furthermore, the mapping
$\M_{\A,x}\ni \theta \rightarrow \theta(U) \in \sigma(U)$
is a continuous map between compact spaces that has dense range.
Since $\theta \in \M_{\A,x}$ is determined by $\theta(U)$,
this map is an injection.
By elementary topology,
the above map is a homeomorphism. Its inverse will be denoted as $\Theta_x$.

Call a map $\Theta$ from a domain $\Omega \subset \bbC$ into $\M_\A$
\textit{pointwise analytic} if $\Theta(z)(a)$ is analytic for $z \in \Omega$
for every $a \in \A$.

\begin{thm}   \label{character}
Let $\X$ be a locally compact Hausdorff space,
$\eta$ a proper continuous map on $\X$ and $x$ a fixed point for $\eta$.
If $\A$ is a topological conjugacy algebra for $(\X, \eta)$, then
there exists a homeomorphism,
\[
\Theta_x: \sigma(U) \longrightarrow \M_{\A, x}
\]
which is pointwise analytic on the interior $\sigma(U)^{\circ}$ of $\sigma(U)$
and satisfies \break $\Theta_x(z)(gU)= g(x)z$ for every $g \in C_0(\X)$.
\end{thm}

\Prf The compact case is known and has been discussed above.
However, the general locally compact case requires an argument.
Indeed, if $\X$ is not compact, then
the corresponding map $\hat{\Theta}_x$ for the canonical unitization $\hat{\A}$ will
do the job provided that we show that every character in $\M_{\A, x}$ comes from one
in $\M_{\hat{\A}, x}$.

To this end, let $\theta \in  \M_{\A, x}$ and let $(e_{\alpha})$ be
the contractive approximate unit of $C_0(\X)$
of all positive functions of norm less than $1$ with compact support.
Let $\lambda= \lim_{\alpha} \theta(e_{\alpha}U)$.
This limit exists because eventually $e_{\alpha}(y) = 1$
on a neighbourhood of $x$, say for $\alpha > \alpha_0$.
Thus if $\alpha, \beta > \alpha_0$, choose $g\in C_0(X)$
so that $g(e_\alpha - e_\beta) = e_\alpha - e_\beta$ and
$g(x)=0$.
Then
\[
 \theta((e_\alpha - e_\beta)U) = \theta(g(e_\alpha - e_\beta)U)
  = \theta(g)  \theta((e_\alpha - e_\beta)U) = 0 .
\]

Note that $|\lambda| \leq \lim \|U^n\|^{1/n}$. Indeed, given $n \in \bbN$,
\begin{eqnarray*}
|\lambda|^n &=&\lim_\alpha |\theta((e_{\alpha}U)^n)|     \\
            &=& \lim_\alpha \Big| \theta\Big( \big(\prod_{m=0}^{n-1}  e_{\alpha}\circ\eta^{(m)} \big)
                       U^n \Big) \Big|  \\
            &\leq& M \| U^n\|,
\end{eqnarray*}
where $M$ is the norm of the embedding of $C(\hat{\X})$ in $\hat{\A}$.
Therefore, $|\lambda|\leq \lim_n \|U^n\|^{1/n}$ and so $\lambda \in \sigma(U)$.

Now notice that for any $g \in C_0(\X)$,
\[
\theta(gU^n) =\lim_{\alpha} \theta(gU^{n-1} e_{\alpha}U)=\theta(gU^{n-1})\lambda, \quad n \in \bbN.
\]
Therefore, $ \theta(gU^n) = g(x)\lambda^n$, $n =1, 2, \dots$, and so $\theta = \theta_{x, \lambda} \in
\M_{\hat{\A}, x}$.

The analyticity of $z \rightarrow \Theta(z)(a)$ on $\sigma(U)^{\circ}$ is
clear for any skew polynomial $a$,
and therefore by approximation, for any $a \in \A$.
\bx

In light of Theorem \ref{character}, given a fixed point $x \in \X$ for $\eta$, the set
\begin{equation} \label{interiorofM}
\{ \theta_{x ,z}\,\mid\, z \in \sigma(U)^{\circ} \}  = \Theta_x(\sigma(U)^{\circ})
\end{equation}
will be simply denoted as $(\M_{\A, x})^{\circ}$.

Let us call a subset of $\M_{\A}$ an \textit{analytic disc} if it is the range
of an injection  $\Phi : (\bbD_s)^{\circ} \to \M_{\A}$, $s>0$,
which is pointwise analytic.
Note that for any
$f \in C_0(\X)$, we have that $\Phi(z)(\overline{f})=\overline{\Phi(z)(f)}$.
So by analyticity, $ \Phi(z)(f)$ must be constant.
Therefore an analytic disc in $\M_{\A}$ is contained in some $\M_{\A, x}$.
Theorem \ref{character} implies that for each fixed
point $x$ of $\eta$, the set $(\M_{\A, x})^{\circ}$ is an a analytic disc
which, by the Open Mapping Theorem, is maximal with that property.

The maximality of $(\M_{\A, x})^{\circ}$ as analytic discs, whenever $x$ is a fixed point for $\eta$, is crucial
for the study of isomorphisms between conjugacy algebras. (The idea of using maximal analytic discs
in the classification of semicrossed products originates in the work of
Hoover \cite{Hoov} and Power \cite{Pow}.)

\section{Isomorphisms between conjugacy algebras}

In order to study isomorphisms between conjugacy algebra we make use of nest representations.

If $\A$
is an algebra, let $\rep_{\fT_2} \A$ will denote the collection of all representations of $\A$ onto $\fT_2$,
the upper triangular $2 \times 2$ matrices.
To each $\pi \in \rep_{\fT_2} \A$ we associate two characters
$\theta_{\pi, 1}$ and $\theta_{\pi, 2}$ which correspond to compressions on the $(1,1)$ and
$(2,2)$-entries, i.e.,
\[
\theta_{\pi, 1}(a)\equiv \langle \pi(a)\xi_i ,  \xi_i\rangle, \quad a \in \A, i=1,2,
\]
where $\{\xi_1,\xi_2\}$ is the canonical basis of $\bbC^2$. If $\gamma:\A_1 \rightarrow \A_2$
is an isomorphism between algebras, then $\gamma$ induces isomorphisms,
\begin{alignat}{2}
\gamma_c&:\M_{\A_1} \to \M_{\A_2}
&\quad\text{by}\quad& \gamma_c(\theta) =  \theta\circ\gamma^{-1} \label{induced1}\\
\gamma_r&:\rep_{\fT_2}  \A_1 \to \rep_{\fT_2}  \A_2
&\quad\text{by}\quad& \gamma_r(\pi) =  \pi\circ\gamma^{-1} \label{induced2},
\end{alignat}
which are compatible in the sense that,
\begin{equation} \label{compat1}
\gamma_c (\theta_{\pi, i}) = \theta_{\gamma_r(\pi),i}, \quad i=1,2,
\end{equation}
for any $\pi \in \rep_{\fT_2}  \A_1 $.

Now assume that $\A$ is a topological conjugacy algebra for $(\X, \eta)$, where $\X$ is
a locally compact Hausdorff space and $\eta$ a proper continuous map. For
$x_1 , x_2 \in \X$, let
\[
\rep_{x_1,x_2}\A \equiv \{\pi \in \rep_{\fT_2} \A : \theta_{\pi, i}
\in \M_{\A, x_i}, i=1,2\} .
\]
Clearly, any element of $\rep_{\fT_2}  \A$ belongs to $\rep_{x,y}\A$ for some
$x,y \in \X$.

\begin{lem} \label{fundrel}
Let  $\X$ be
a locally compact Hausdorff space, $\eta$ a proper continuous map on $\X$ and
$\A$ a topological conjugacy algebra for $(\X, \eta)$. Assume that $x,y \in \X$
are not fixed points for $\eta$ and let $\pi \in \rep_{x,y}\A$. Then, $y=\eta(x)$.
\end{lem}

\Prf By assumption, $\theta_{\pi, 1}= \theta_{x,0}$ and $\theta_{\pi, 2}= \theta_{y,0}$
and so $\theta_{\pi, 1}(gU) = \theta_{\pi, 2}(gU)=0$, for any $g \in C_0(\X)$.
Therefore for each $g \in C_0(\X)$ there exists $c_g \in \bbC$
so that
\[
\pi(gU)=
\left(\begin{array}{cc} 0 & c_g \\  0 & 0 \end{array} \right).
\]
By Proposition \ref{Cohen} there exists at least one $g \in C_0(\X)$ so that
$c_g\neq0$, or otherwise the range of $\pi$ would be commutative.
Applying $\pi$ to $gUf= (f\circ\eta) gU$
for $f \in C_0(\X)$ and this particular $g$, we get
\begin{equation*}
\left(\begin{array}{cc} 0 & c_g \\ 0 & 0 \end{array} \right)
\left(\begin{array}{cc} f(x) & t \\ 0 & f(y)\end{array} \right)
=
\left(\begin{array}{cc} f(\eta(x)) & t' \\ 0 & f(\eta(y)) \end{array} \right)
\left(\begin{array}{cc} 0 & c_g \\  0 & 0 \end{array} \right)
\end{equation*}
for some $t, t' \in \bbC$, depending on $f$.
By comparing $(1,2)$-entries, we obtain,
\[
f(y) = f (\eta(x)) \qforal f \in C(\X),
\]
i.e., $y = \eta(x)$, as desired.
\bx

\begin{eg} \label{example1}
$\rep_{x,\eta(x)}\A$ is not empty when $\eta(x) \ne x$.
Let
\[
\rho(f)=
\left(\begin{array}{cc} f(x) & 0 \\ 0 & f(\eta(x))\end{array} \right), \quad
\rho(fU)=
\left(\begin{array}{cc} 0 & f(x) \\ 0 & 0 \end{array} \right)
\]
and $\rho(fU^n)=0$, $n\geq2$, $f \in C_0(\X)$. Extend $\rho$ by linearity to a
map $\rho:P^{\infty}(X, \eta)\rightarrow\fT_2$. Then, $\rho\circ\Delta
\in \rep_{x,\eta(x)}\A$, where $\Delta$ is the Fourier series defined by (\ref{Delta}).
\end{eg}

We do not know if the analog of Lemma \ref{fundrel} is valid
when $y$ is a fixed point without assuming the continuity of
the representation $\pi$. In order to overcome this difficulty, we introduce
a global object.

\begin{defn}  \label{pencil}
Let $\X$ be a locally compact Hausdorff space, $\eta$ a proper continuous map on $\X$ and
$\A$ a topological conjugacy algebra for $(\X, \eta)$. Assume that $x,y \in \X$
so that $\eta(y) = y$ but $\eta(x)\neq x$. A \textit{pencil
of nest representations} for $\A$ is a set $\P_{x,y} \subseteq \rep_{x,y}\A$ which satisfies
\[
 \{\theta_{\pi,2} : \pi \in \P_{x,y} \}= (\M_{\A, y})^{\circ},
\]
where $(\M_{\A, y})^{\circ}$ is defined by (\ref{interiorofM}).
\end{defn}

\begin{lem} \label{fundrelp}
Let $\X$, $\eta$, $x,y \in \X$ and $\A$ be as in Definition \ref{pencil}
and let $\P_{x,y}$ be a pencil of representations for $\A$ . Then, $y =\eta(x)$.
\end{lem}

\Prf Since $\P_{x,y}$ is a pencil, there exists $\pi \in \P_{x,y}$ so that
$\theta_{\pi, 1}= \theta_{x,0}$ and $\theta_{\pi, 2}= \theta_{y,0}$. The rest of the proof
now is identical to that of Lemma~\ref{fundrel}.
\bx

\begin{eg} \label{expen}
Assume that $\A$ is a conjugacy algebra for a dynamical system $(\X, \eta)$
and $x \in \X$ satisfies $x\neq\eta(x)$ and $\eta^{(2)}(x)=\eta(x)$. We show that there exists
a pencil of representations for $\A$ of the form $\P_{x , \eta(x)}$. \smallbreak

First note that it suffices to consider only compact spaces $\X$.
Recall that $\liminf_n\|E_n\|^{-1/n}\geq r$, where $r$ is the spectral radius
of $U$ (and so $\sigma(U)=\bbD_r$).
If $|z| < r$, then we define,
\[
 \pi_z (f)=
 \left(\begin{array}{cc} f(x) & 0 \\ 0 & f(\eta(x))\end{array} \right), \quad \mbox \quad
 \pi_z (U)=
 \left(\begin{array}{cc} 0 & z \\ 0 & z\end{array} \right).
\]
For an arbitrary $a \sim \sum_n  E_n (a)U^n$, define
\begin{eqnarray*}
 \pi_z(a) &=& \sum_n \pi_z(E_n(a)) \pi_z(U)^n    \\
 &=&
 \left(\begin{array}{cc} E_0(a)(x) & \sum_{n\geq1} E_n(a)(x)z^n \\ 0 &
 \sum_{n\geq0} E_n(a)(\eta(x))z^n \end{array} \right).
\end{eqnarray*}
Since $|z|<r$, $\pi_z (a)$ is well defined for any $a \in \A$. Furthermore,
it is easy to see that $\pi_z(U^n)=\pi_z(U)^n$, $n \in \bbN$, and also
$\pi_z(f\circ \eta)\pi_z(U) = \pi_z(U)\pi_z(f)$. From this, it easily follows that $\pi_z$
is an algebra homomorphism that maps onto $\fT_2$.
The desired pencil of representations is therefore
\[
 \P_{x, \eta(x)}\equiv \{ \pi_z  :  z \in \sigma(U)^{\circ}\} ,
\]
where $\pi_z$ are as above. \smallbreak

More generally, suppose that $\pi \in \rep_{x,y}\A$ is continuous
where $y = \eta(x) = \eta(y)$ is a fixed point and $\theta_{\pi,2} = \theta_{y,z}$
for some $z \in \bbD_r = \sigma(U)^\circ$.
Then there is a scalar $\alpha$ so that
\[
 \pi(U) = \begin{pmatrix} \theta_{x,0}(U) & \alpha \\ 0 & \theta_{y,z}(U) \end{pmatrix}
 =  \begin{pmatrix}0 & \alpha \\ 0 & z \end{pmatrix}.
\]
Thus
\[
 \pi(f U^k) =  \begin{pmatrix}f(x) & 0 \\ 0 & f(y) \end{pmatrix}
 \begin{pmatrix} 0 & \alpha z^{k-1} \\ 0 & z^k \end{pmatrix}
 =  \begin{pmatrix}0 & \alpha f(x) z^{k-1} \\ 0 & f(y) z^k \end{pmatrix}.
\]
By continuity and the fact that $|z| < r$, we obtain that
\[
 \pi(a) =  \begin{pmatrix} E_0(a)(x) & \alpha \sum_{k=1}^\infty E_k(a)(x) z^{k-1} \\
 0 & \sum_{k=0}^\infty E_k(a)(y) z^k  \end{pmatrix}.
\]

For $|z|=r$, we do not have a complete picture.

When $\pi$ is discontinuous, it is possible that
$\pi(U) =  \left(\begin{smallmatrix} 0&0\\0&z \end{smallmatrix}\right)$
yet there is some element $a \in \A$ such that
$\pi(A) =  \left(\begin{smallmatrix} 0&1\\0&0 \end{smallmatrix}\right)$.
We do not know whether such representations actually exist.
Also when $\eta(x)=x$ is a fixed point, we do not know if
$\rep_{x,x} \A$ is non-empty.  Such a representation
is necessarily discontinuous since functions are sent to the
scalars, and so the skew polynomials are sent onto
the abelian algebra generated by $\pi(U)$.
\end{eg}

We are ready to state and prove the main result of this section.

\begin{thm}   \label{main}
Let $\X_i$ be a locally compact Hausdorff space and let
$\eta_i$ a proper continuous map on $\X_i$, for $i=1,2$. Then the dynamical
systems $(\X_1, \eta_1)$ and $(\X_2, \eta_2)$ are conjugate
if and only if some topological conjugacy algebra for
$(\X_1, \eta_1)$ is isomorphic as an algebra to some
topological conjugacy algebra for
$(\X_2, \eta_2)$.
\end{thm}

\Prf If the two systems are conjugate, then the algebras of
Example~\ref{rightone} are easily seen to be isomorphic.

Conversely, suppose that $\A_i$ are conjugacy algebras for $(\X_i, \eta_i)$, $i=1,2.$, and
that there exists an algebra homomorphism $\gamma: \A_1\rightarrow\A_2$.
Then $\gamma$ induces a homeomorphism $\gamma_c$ of
$\M_{\A_1}$ onto $ \M_{A_2}$ by $\gamma_c(\theta) = \theta\circ \gamma^{-1}$.
It is elementary to verify that
$\gamma_c$ preserves analytic discs and therefore establishes a bijection between the maximal
analytic discs of $\A_1$ and $\A_2$. This bijection extends to a bijection between
their closures and therefore to a bijection between the
collections $\{\M_{\A_1 , x} : x \in \X_1 \}$ and
$\{\M_{\A_2 , y} : y \in \X_2 \}$.
In other words, for each $x \in \X_1$ there exists
a $\gamma_s(x) \in \X_2$ so that
\begin{equation} \label{compat2}
\gamma_c(\M_{\A_1 , x})=\M_{\A_2 , \gamma_s(x)}.
\end{equation}
We have therefore defined a map $\gamma_s: \X_1 \rightarrow \X_2$, which maps fixed points to fixed points and
satisfies
\begin{equation} \label{dif}
f(\gamma_s(x))=(\theta_{x,0}\circ \gamma^{-1})(f)
\end{equation}
for all $x \in \X_1$ and $f \in C_0(\X_2)$.
Notice that if $(x_i)_i$ is a net in converging to some $x \in \X_1$, then
(\ref{dif}) shows that $\big(f(\gamma_s(x_i))\big)_i$ converges to $f(\gamma_s(x))$, for all $f \in
C_0(\X_2)$, and so $(\gamma_s(x_i))_i$ converges to $\gamma_s(x)$. Hence, $\gamma_s $ is continuous.
Repeating the above arguments
with $\gamma^{-1}: \A_2\rightarrow\A_1$ in the place of $\gamma$, we obtain
that $\gamma_s: \X_1 \rightarrow \X_2$ has a continuous inverse and is therefore a homeomorphism.
Furthermore,$\gamma_s$ maps
the fixed point set of $\eta_1$ onto the fixed point set of $\eta_2$. Finally,
\smallbreak
\noindent
\textbf{Claim 1:}\, If $x \in \X_1$ is not a fixed point for $\eta_1$, then
\[
\gamma_r(\rep_{x, \eta_1(x)}\A_1)
\subseteq \rep_{\gamma_s(x) , \gamma_s(\eta_1(x))} \A_2 .
\]
Pick a representation $\pi \in \rep_{x, \eta_1(x)}\A_1$.
By (\ref{compat1}), we have $\theta_{\gamma_r(\pi),1} = \gamma_c(\theta_{\pi, 1})$
and so $\theta_{\gamma_r(\pi),1}=\gamma_c(\theta_{x,0})$. By (\ref{compat2}),
$\gamma_c(\theta_{x,0}) \in \M_{\A_2, \gamma_s(x)}$ and so $\theta_{\gamma_r(\pi),1}
\in \M_{\A_2, \gamma_s(x)}$. A similar argument shows that $\theta_{\gamma_r(\pi),2}
\in \M_{\A_2, \gamma_s(\eta_1(x))}$ and this proves the claim.\vspace{2mm}

We now show that
$\gamma_s$ implements the desired conjugacy between $(\X_1, \eta_1)$ and $(\X_2, \eta_2)$,
i.e.,
\begin{equation} \label{funeq}
\gamma_s(\eta_1(x))=\eta_2(\gamma_s(x)), \, \mbox{for all } x \in \X_1.
\end{equation}

Since $\gamma_s$ maps fixed points to fixed point, verifying (\ref{funeq}) becomes trivial
in that case. We therefore pick an $x \in \X$ with $\eta_1(x)\neq x$
and we examine two cases.

For the first case assume that $\eta_1^{(2)}(x)\neq \eta_1(x)$. In that case,
pick a representation $\pi \in \rep_{x ,\eta_1(x)} \A_1$.
Combining Claim 1 with Lemma \ref{fundrel}, we obtain that
$\eta_2( \gamma_s(x) )= \gamma_s(\eta_1(x))$, which proves (\ref{funeq}) in the first
case.

For the second case assume that $\eta_1^{(2)}(x)= \eta_1(x)$ and let
$\P_{x , \eta_1(x)}$ be a pencil of representations for $\A_1$ as in Example \ref{expen}.
By Claim 1,
\[
\gamma_r(\P_{x, \eta_1(x)}) \subseteq \rep_{\gamma_s(x) , \gamma_s(\eta_1(x))} \A_2.
\]
Since $\gamma_c$ preserves maximal analytic discs, $\gamma_r$ preserves pencils of representations and
so $\gamma_r(\P_{x, \eta_1(x)}) $ is a pencil of the
form $\P_{\gamma_s(x) , \gamma_s(\eta_1(x))} $. By Lemma \ref{fundrelp},
we have $\eta_2( \gamma_s(x) )= \gamma_s(\eta_1(x))$, which proves (\ref{funeq}) in the
last remaining case. This proves the Theorem.
\bx

\begin{cor}
Let $\X_i$ be a locally compact Hausdorff space and let
  $\eta_i$ a proper continuous map on $\X_i$, for $i=1,2$. Then the dynamical
  systems $(\X_1, \eta_1)$ and $(\X_2, \eta_2)$ are conjugate
  if and only if the semicrossed products $C_0(\X_1) \times_{\eta_1}\bbZ^+$
  and $C_0(\X_2) \times_{\eta_2}\bbZ^+$ are isomorphic as algebras.
  \end{cor}

%%%%%%%%%%%%%%%%%%%%%%%%%%%%%%%%%%%%%%%%%%%%%
 
 \section{Isomorphisms between subalgebras of conjugacy algebras}

The techniques of the previous section are widely applicable and suggest new
avenues of investigation in the classification of non-selfadjoint algebras.
In this section we explore just one such possible direction.

Let $\X$ be a compact Hausdorff space, and let $\eta$ be a continuous map on $\X$.
Let $\S \subseteq C(\X)$ be a uniform algebra, and assume that $\eta$ leaves
$\S$ invariant, i.e., if $f \circ \eta \in \S$ for every $f \in \S$. If $\A$ is
a conjugacy algebra for $(\X, \eta)$, then $\A(\S)$ will denote the norm closed
subalgebra of $\A$ generated by all polynomials of the form
$\sum _{n}  f_n U^n$ where $f_n \in \S$ for all $n\ge0$. Such a subalgebra
inherits some of the properties of the enveloping conjugacy algebra. For instance,
the formal power series expansion of (\ref{powerseries}) is always valid for $\A(\S)$.
If the characters of the uniform algebra $\S$ coincide
with the point evaluations on $\X$, then
Theorem \ref{character} is also valid for $\A(\S)$.

For simplicity, we restrict our attention to the class $\K$ of subsets
$\bbK$ of $\bbC$ of the form $\bbK = \overline{G}$, where $G$
is a Cauchy domain (i.e. a connected open set in $\bbC$ with
boundary equal to the union of finitely many disjoint Jordan curves).
Let $A(\bbK)$ denote the  algebra of functions which are
continuous on $\bbK$ and analytic on its interior.
Clearly one can apply our methods to more complicated regions,
but the principle is clear for this class.

It turns out that for our problem, the only pathology occurs when
$G$ is simply connected, which reduces to the case of the unit disk.

\begin{thm} \label{analytic}
Let $\bbK_i \in \K$,
let $\eta_i :\bbK_i \rightarrow \bbK_i$ be a continuous map
which is analytic on the interior, and
let $\A_i$ be a conjugacy algebra for $(\bbK_i, \eta_i)$ for $i = 1,2$.
If the algebras $\A_1(A(\bbK_1))$ and $\A_2(A(\bbK_2))$
are isomorphic, then either
\begin{enumerate}
\item $\eta_1$ and $\eta_2$ are analytically conjugate;
i.e. there exists a homeomorphism $\gamma_s: \bbK_1 \to \bbK_2$ which is
analytic  in the interior and satisfies
$\gamma_s\circ \eta_1 = \eta_2\circ\gamma_s$; or

\item $\bbK_i$ are simply connected, $\eta_i$ are homeomorphisms
with a unique fixed point in the interior, and $\eta_2$ is
analytically conjugate to $\eta_1^{-1}$.
\end{enumerate}
\end{thm}

In the rest of this section, we develop the tools we need to
prove this theorem.

\begin{eg}
Let $\bbK \in \K$ and let $\eta: \bbK \rightarrow \bbK$ be a continuous map
which is analytic on the interior.
Let $\H_x = l^2(\bbN)$ be the space of square summable sequences
$\xi = (\xi_n)_{n=0}^{\infty}$
and for each $x \in \bbK$ define
\[
\pi_x(f)\xi = (f(x)\xi_0,(f\circ \eta)(x)\xi_1, (f\circ\eta^{(2)})(x)\xi_2, \dots)
\qfor f \in A(\bbK),
\]
and
\[
V_x\xi = (\xi_1, \xi_2, \xi_3, \dots).
\]
The norm closed operator
algebra acting on $\oplus_{x \in \bbK}\H_x$ and generated by the
operators,
\begin{equation*}
\oplus_{x \in \bbK} \pi_x(f), \quad
 \oplus_{x \in \bbK} \pi_x(g)V_x \qfor f,g \in A(\bbK),
\end{equation*}
is an algebra of the form $\A_{\bbK, \eta}(A(\bbK))$, where $\A_{\bbK, \eta}$
is as in Example \ref{rightone}.
\end{eg}

\begin{eg} \label{finalex}
Let $\bbK$, $\H_x$ and $\pi_x(f)$ be as above, and let
$U_x$ be the forward shift
\[
U_x\xi = (0, \xi_0, \xi_1, \xi_2, \dots).
\]
We (temporarily) denote as $\alg(A(\bbK), \eta)$ the norm closed operator
algebra acting on $\oplus_{x \in \bbK}\H_x$ and generated by the
operators,
\begin{equation}
\oplus_{x \in \bbK} \pi_x(f), \quad
 \oplus_{x \in \bbK} U_x\pi_x(g) \qfor f,g \in A(\bbK).
\end{equation}
As we explained in Example \ref{exsemicrossed}, our theory is applicable to the opposite algebras
$(\alg(A(\bbK), \eta)_{op}$ but can be used to classify the algebras
$\alg(A(\bbK), \eta)$ as well.
\end{eg}

We need to make some clarifications regarding the example above and the theory in \cite{BP}.
Let $\bbD \subseteq \bbC$ denote the closed unit disk, $\bbT$ its boundary and let
$\eta : \bbD \rightarrow \bbD$ be a \textit{homeomorphism} which is analytic on the interior.
Let $S$ be the forward unilateral shift on $\oplus_{n=0}^{\infty} L^{2}(\bbT)$, and for
$f \in A(\bbD)$, let $\phi(f)$ denote the operator on $\oplus_{n=0}^{\infty} L^{2}(\bbT)$
defined by
\[
  \phi(f)(\xi_0, \xi_1, \xi_2 , \dots)=
  (f\xi_0, (f\circ \eta)\xi_1, (f\circ\eta^{(2)})\xi_2, \dots).
\]
In \cite{BP}, the norm closed operator algebra generated by $S$ and the operators
$\phi(f)$ for $f \in \A(\bbD)$ is denoted as $\fA_{\eta}$.

\begin{prop}  \label{clarify}
Let $\eta : \bbD  \to \bbD$ be a homeomorphism which is analytic on the interior.
Then the algebras $\fA_{\eta}$ and $\alg(A(\bbD), \eta)$ defined above are isomorphic.
\end{prop}

\Prf Let $P(A(\bbD), \eta)$ be the skew polynomial algebra with coefficients in $A(\bbD)$.
The algebra $P(A(\bbD), \eta)$ can be equiped with three norms
which are defined as
\begin{align*}
 \|p\|_1&= \Big\| \sum_i \big(\oplus_{x \in \bbT} U_x\big)^i \big(\oplus_{x \in \bbT}\pi_x(f_i)\big) \Big\| \\
 \| p \|_2 &= \Big\| \sum_i \big(\oplus_{x \in \bbD} U_x\big)^i \big(\oplus_{x \in \bbD}\pi_x(f_i)\big) \Big\|
\end{align*}
and
\[
  \|p\|_3= \| \sum_i  S^i \phi(f_i) \|,
\]
for $p\equiv \sum_i U^if_i \in P(A(\bbD), \eta) $.
By \cite[Proposition IV.1]{BP}, we have
\[
  \|p\|_1 \leq \|p\|_2 \leq \|p\|_3
\]
for any $p \in P(A(\bbD), \eta)$. Clearly the result will follow
if we show that the above inequalities are actually equalities.

Towards this end, note that \cite[Proposition II.7]{Pet} implies that \break
$(P(A(\bbD), \eta),  \| \,.\|_1)$ embeds isometrically into the crossed product C*-algebra
$C(\bbT) \times_{\eta} \bbZ$.
Therefore, if $(W_1, W_2) $ is a pair of unitary operators which is universal
with  the property that $W_1W_2 = W_2 \eta(W_1)$, then
\[
  \|p\|_1= \| \sum_i  W_2^i f_i(W_1) \|,
\]
for any
$p \in P(A(\bbD), \eta)$. But then, \cite[Theorem II.4]{BP} implies that
$\| p\|_3  \leq\|p\|_1$, $p \in P(A(\bbD), \eta)$, and the conclusion follows.
\bx

In light of Proposition \ref{clarify} and the theory in \cite{BP},
we will no longer use the notation $\alg( A(\bbK), \eta)$ of
Example \ref{finalex}, and for the rest of the paper we will refer to
this algebra as $A(\bbK) \times_{\eta} \bbZ^{+}$.

\begin{defn}
Let $\A$ be a Banach algebra with character space $\M_{\A}$. We say
that an open subset $\Omega\subseteq \M_{\A}$ is an \textit{analytic set} if there exists
a pointwise analytic injective map
$\Theta:G \to \M_{\A}$ from a domain  $ G \subset \bbC$
with range $\Omega$.
\end{defn}

Assume that $\bbK \in \K$, $\eta: \bbK \rightarrow \bbK$, $\eta \neq id$, is a continuous map
which is analytic on the interior, and let $\A$ be a conjugacy algebra for $(\bbK, \eta)$.
The maximal analytic sets in the character space of $\A$ coincide with the maximal analytic
discs and were identified in the discusion following Theorem \ref{character}. In the case
of $\A(A(\bbK))$ we have one additional analytic set.

\begin{prop}   \label{maxdiscs}
Let $\bbK \in \K$, let $\eta: \bbK \rightarrow \bbK$ be a continuous map
 which is analytic on the interior, and is not the identity map.
Let $\A$ be a conjugacy algebra for $(\bbK, \eta)$.
Then the maximal analytic sets
in the character space of  $\A(A(\bbK))$
are the sets $(\M_{\A(A(\bbK)), x})^{\circ}:=\{ \theta_{x,z} : |z| < r \}$
for fixed points $\eta(x) = x$,
where $r$ is the spectral radius of $U$, and
$(\M_{\A(A(\bbK))}^0)^{\circ} := \{ \theta_{x,0}  : x \in \bbK^{\circ} \}$.
\end{prop}

\Prf Suppose that $G$ is a domain in $\bbC$ and
\[
  \Theta : G \to \M_{\A(A(\bbK))}
\]
is pointwise analytic.
We will show that either
the range of $\Theta$ is contained in one of the sets
$(\M_{\A(A(\bbK)), x})^{\circ}$, $x \in \bbK$,
or otherwise, $\Theta(\bbK) \subseteq (\M_{\A(A(\bbK))}^0)^{\circ}$.
This will clearly imply the result.

Let us write $\Theta(z) = \theta_{\delta(z) ,\iota(z)}$ for $z \in G$,
where $\delta$ and $\iota$ map $G$ into $\bbK$ and $\bbD_r$ respectively.
Note that the map $z \to \Theta(z)(f)= f(\delta(z))$ for $z \in G$ is
analytic. By considering the identity function $f(x)=x$, we conclude that
$\delta$ is analytic.
We distinguish two cases.

Assume first that there exists $z_0 \in G$
so that $\Theta(z_0)(U)\neq0$. By continuity,
there exists $\eps >0$ so that $\Theta(z)(U)\neq0$,
for all $z \in B_\eps(z_0)$.
Hence $\iota(z) \neq 0$ and so $\delta (z)$ is a fixed point for $\eta$
for all $z \in B_\eps(z_0)$.
If $\delta$ were not constant, $\delta(B_\eps(z_0))$ would be
an open set of fixed points of $\eta$.
By analyticity, $\eta = \id$ contrary to hypothesis.
Hence $\delta$ is constant.

For the second case, there is no $z \in G$
so that $\Theta(z)(U)\neq0$,
and so $\Theta(G) \subseteq (\M_{\A(A(\bbK))}^0)^{\circ}$.
\bx

For the rest of the section, the closures of the maximal analytic discs
$(\M_{\A(A(\bbK)), x})^{\circ}$, $x \in \bbK$, and $(\M_{\A(A(\bbK))}^0)^{\circ}$
will be denoted as
$\M_{\A(A(\bbK)), x}$, $x \in \bbK$, and $\M_{\A(A(\bbK))}^0$
respectively.

\begin{thm} \label{generalregion}
Let $\bbK_i \in \K$,
let $\eta_i :\bbK_i \rightarrow \bbK_i$ be a continuous map
which is analytic on the interior and
let $\A_i$ be a conjugacy algebra for $(\bbK_i, \eta)$, $i = 1,2$.
Assume furthermore that $\bbK_1$ is not simply connected.
If the algebras $\A_1(A(\bbK_1))$ and $\A_2(A(\bbK_2))$
are isomorphic, then there exists a homeomorphism
$\gamma_s: \bbK_1 \rightarrow \bbK_2$ which is
analytic  in the interior and satisfies
$\gamma_s\circ \eta_1 = \eta_2\circ\gamma_s$.
\end{thm}

\Prf Assume that there exists an algebraic isomorphism
$\gamma$ from $\A_1(A(\bbK_1))$ to $\A_2(A(\bbK_2))$; and let
$\gamma_c: \M_{\A_1(A(\bbK_1)) } \rightarrow \M_{\A_2(A(\bbK_2)) }$
be as in (\ref{induced1}).
In this and other special cases of Theorem~\ref{analytic}(i),
the goal is to show that $\gamma_c$ is a homeomorphism of
$\M_{\A_1(A(\bbK_1))}^0$ onto $\M_{\A_2(A(\bbK_2))}^0$.
Then one proceeds as in Theorem \ref{main}.

Assume that $\eta_1$ is not the identity map.
Evidently $\gamma_c$ maps maximal analytic sets to maximal analytic sets.
The maximal analytic sets in $\M_{\A_1(A(\bbK_1)) }$ consist of
open disks over each fixed point and one set homeomorphic to $\bbK_1^\circ$
which is not simply connected.
Hence $\bbK_2$ also has such a set, which is the image under $\gamma_c$.
By continuity, $\gamma_c$ maps the closure $\M_{\A_1(A(\bbK_1))}^0$
onto $\M_{\A_2(A(\bbK_2))}^0$.
Hence there is a bijection between the collections
$\{\M_{\A_1(A(\bbK_1)),x} : x \in \bbK_1 \}$, and
$\{\M_{\A_2(A(\bbK_2)) , x} : x \in \bbK_2 \}$, which allows us to define a map
$\gamma_s :\bbK_1\to \bbK_2$ so that
$\gamma_c(\M_{\A_1(A(\bbK_1)) , x})=\M_{\A_2(A(\bbK_2)) , \gamma_s(x)}$.
One now proceeds as in the proof of Theorem \ref{main} to show that
$\gamma_c$ is a homeomorphism that intertwines $\eta_1$ and $\eta_2$.
By equation (\ref{dif}) applied to the identity function $\bz(x)=x$,
we obtain
\begin{equation} \label{difa}
\gamma_s(x) = \theta_{x,0}(\gamma^{-1} \bz) = g(x)
\end{equation}
where $g = E_0 (\gamma^{-1} \bz)$ belongs to $A(\bbK_1)$.
Hence $\gamma_c$ is analytic on $\bbK_1^\circ$.

It remains to deal with the case where $\eta_1$ is the identity map.
Then $\A(\bbK_1) \times_\id \A(\bbD) \simeq A(\bbK_1 \times \bbD)$
is abelian.  Hence $\A(\bbK_2) \times_{\eta_2} \A(\bbD)$
is also abelian.  Therefore $\eta_2$ is the identity map as well.
Any algebra isomorphism from $A(\bbK_1 \times \bbD)$ onto
$A(\bbK_2 \times \bbD)$ is determined by a homeomorphism
of the maximal ideal spaces $\bbK_i \times \bbD$ which is
biholomorphic on the interior.
By \cite[Prop. 2]{Lig} such a map must be a product map.
When $\bbK_i$ are simply connected, they are clearly biholomorphic;
while if there are not simply connected, the product
map must carry $\bbK_1^\circ$ onto $\bbK_2^\circ$.
This map extends to be continuous on the boundary,
establishing the desired equivalence.
\bx

Any simply connected region $\bbK$ with a Jordan curve as boundary is
conformally equivalent to the disk $\bbD$; i.e. there is a homeomorphism
of $\bbK$ onto $\bbD$ which is analytic on the interior.
So we can consider only the region $\bbD$ for the remainder.
A continuous map $\eta:\bbD \to \bbD$ which is analytic on the interior
is called \textit{elliptic} if $\eta$ has exactly one fixed point which is located
in the interior of $\bbD$.  Indeed, if $\eta \ne \id$, then a fixed point interior
to $\bbD$ implies the uniqueness of the fixed point by Schwarz's Lemma.

 \begin{thm}    \label{secmain}
Let $\eta_i :\bbD \rightarrow \bbD$ be a continuous map which is analytic on the interior, and
let $\A_i$ be a conjugacy algebra for $(\bbD, \eta)$, $i = 1,2$. Assume further that
$\eta_1$ is not elliptic. If the algebras $\A_1(A(\bbD))$ and $\A_2(A(\bbD))$
are isomorphic, then there exists a homeomorphism \break $\gamma_s: \bbD \rightarrow \bbD$ which is
analytic on the interior and satisfies
$\gamma_s\circ \eta_1 = \eta_2\circ\gamma_s$.
\end{thm}

\Prf  Assume that there exists an algebraic isomorphism \break
$\gamma: \A_1(A(\bbD)) \rightarrow   \A_2(A(\bbD))$, and let
$\gamma_c: \M_{\A_1(A(\bbD)) } \rightarrow \M_{\A_2(A(\bbD)) }$ be as in
(\ref{induced1}).
Again we will show that $\gamma_c( \M_{\A_1(A(\bbD))}^0) = \M_{\A_2(A(\bbD))}^0$.
As $\eta_1$ is not elliptic, there exists a fixed point for $\eta_1$ on the
boundary of $\bbD$, say $x_0$. This implies that $(\M_{\A_j(A(\bbD))}^0)^{\circ}$
is the only maximal analytic disc, whose boundary has non-empty
intersection with some other analytic disc (namely $\{ \theta_{x_0,z} : |z|<r\}$).
This property will be preserved by $\gamma_c$. Hence, $\eta_2$ is not elliptic
and $\gamma_c( \M_{\A_1(A(\bbD))}^0) = \M_{\A_2(A(\bbD))}^0$.
The analyticity of $\gamma_c$ follows as in the proof of Theorem~\ref{generalregion}.
\bx

\begin{cor} \label{maindisc1}
Let $\eta_i: \bbD \rightarrow \bbD$ be continuous map which are analytic on the interior, $i=1,2$.
Assume further that $\eta_1$ is not elliptic.
The semicrossed products $A(\bbD) \times_{\eta_1}\bbZ^{+}$ and
$A(\bbD) \times_{\eta_2}\bbZ^{+}$ are isomorphic as algebras if and only if
there exists a homeomorphism $\gamma_s: \bbD \rightarrow \bbD$ which is analytic
in the interior and satisfies
$\gamma_s\circ \eta_1 = \eta_2\circ\gamma_s$.
\end{cor}

To complete the classification, we have to examine the elliptic case.
To deal with this case, we need to be able to reduce to a situation
algebraic isomorphisms are continuous.
The key is to show that the ideal $\ker \Delta$ is determined
algebraically by the characters.

It is easy to see that if $\A$ is a conjugacy algebra for $(\bbK, \eta)$, where $\eta$ is analytic
in the interior of $\bbK$, then $\ker \Delta = \bigcap_{n\ge1} \J^n$ where
$\J = \ker E_0 = \bigcap_{x \in \bbK^\circ} \ker \theta_{x,0}$.
Except when $\bbK$ is simply connected and $\eta$ has a unique
fixed point $x_0$ and $x_0 \in \bbK^\circ$, the maximal analytic
set $(\M_{\A(A(\bbK))}^{0})^{\circ} = \{ \theta_{x,0} : x \in \bbK^\circ \}$
is distinguished from all others.
In the elliptic case, there are two maximal analytic disks which intersect
each other, namely $(\M_{\A(A(\bbK))}^{0})^{\circ}$ and
$(\M_{\A(A(\bbK)),x_0})^{\circ} = \{ \theta_{x_0,z} : |z| < r \}$.
Let $\J$ and $\K$ be the ideals
obtained from the intersection of their kernels.

\begin{prop} \label{kerDelta}
Let $\eta: \bbD \rightarrow \bbD$ be a continuous map which is analytic on the interior,
and assume further that $\eta_1$ is elliptic.
Let $\A$ be a conjugacy algebra for $(\bbD,\eta)$. Then
\[
 \ker \Delta = \bigcap_{n\ge1} \J^n =
 \bigcap_{n\ge1} \J^n + \bigcap_{n\ge1} \K^n ,
\]
where $\J$ and $\K$ are as above.
\end{prop}

\Prf  Assume that there is a unique fixed point $x_0$ for $\eta$
which is inside the disk.
Then there are exactly two maximal analytic disks.
They intersect in interior points, and are not distinguished
topologically.
Consider the ideal $\K$ defined above.
It is clear that
\[ \K = \{ a : \Delta(a) = \sum_k f_k U^k \AND f_k(x_0)=0 \FOR k \ge 0 \} . \]
We claim that
\[
 \K^n \subset \{ a :
 f_k \text{ has a zero of order }n \text{ at }x_0 \FOR k \ge 0 \}
 =: \K_n.
\]
It suffices to show that $\K_n \K \subset \K_{n+1}$.

Write $\Delta(a) = \sum_k f_k U^k$ and $\Delta(b) = \sum_k g_k U^k$,
where $a \in \K_n$ and $b \in \K$.  Then
\[ \Delta(ab) = \sum_{n\ge0} \sum_{k=0}^n f_k(g_{n-k}\circ \eta^{(k)}) U^n .\]
Since $\eta(x_0)=x_0$, $g_{n-k}\circ \eta^{(k)}$ has a zero at $x_0$
and each $f_k$ has a zero of order at least $n$.
So the product has a zero of order $n+1$ for every term in this sum.

It now follows that
\[
 \bigcap_{n\ge1} \K^n \subset \bigcap_{n\ge1} \K_n = \ker \Delta.
\]
Any isomorphism carries maximal analytic sets to maximal analytic sets.
Whether it fixes them or switches them, one obtains
$\ker \Delta$ algebraically by forming two ideals $\J$ and $\K$
from the intersection of their kernels, and obtaining
\[
 \ker \Delta = \bigcap_{n\ge1} \J^n + \bigcap_{n\ge1} \K^n .
\]
\upbx

The crucial application is the following evident consequence.

\begin{cor}\label{kerDeltacor}
Let $\bbK_i \in \K$,
let $\eta_i :\bbK_i \rightarrow \bbK_i$ be a continuous map
which is analytic on the interior, and
let $\A_i$ be a conjugacy algebra for $(\bbK_i, \eta_i)$, $i = 1,2$.
If the algebras $\A_1(A(\bbK_1))$ and $\A_2(A(\bbK_2))$
are isomorphic, then the conjugacy algebras
$\A_1(A(\bbK_1))/ \ker \Delta_1$ and $\A_2(A(\bbK_2))/ \ker \Delta_2$
are isomorphic.
\end{cor}

Now we will establish an automatic continuity result in the case
in which $\Delta$ is injective.
Recall that if $\phi : \A \to \B$ is an epimorphism between
Banach algebras, then the \textit{separating space} of $\phi$ is the
two-sided closed ideal of  $\B$ defined as
\[
\S (\phi ) \equiv \big\{ b \in \B  : \exists \{a_n\}_n \subseteq \A
\mbox{ such that } a_n \rightarrow 0 \mbox{ and } \phi(a_n )\rightarrow b\big\}.
\]
Clearly the graph of $\phi$ is closed if and only if $\S (\phi )= \{0\}$.
Thus by the closed graph theorem, $\phi$ is continuous
if and only if $\S (\phi ) = \{0\}$.

The following is an adaption of \cite[Lemma 2.1]{Sinclair} and was used in
\cite{DHK} for the study of isomorphisms between limit algebras.

\begin{lem}[Sinclair]     \label{Sinclair}
Let  $\phi : \A \to \B$ be an epimorphism between
Banach algebras and let $ \{B_n\}_{n\in \bbN}$
be any sequence in $\B$.
Then there exists $n_0 \in \bbN$ so that for all $n \geq n_0$,
\[
\overline{ B_1 B_2 \dots B_n \S (\phi ) } = \overline{ B_1 B_2 \dots B_{n+1}  \S (\phi ) }
\]
and
\[
\overline{ \S (\phi ) B_{n} B_{n-1} \dots B_1 } = \overline{\S (\phi ) B_{n+1} B_n \dots B_1 }.
\]
\end{lem}

\begin{cor} \label{auto_cont}
Let $\bbK_i \in \K$,
let $\eta_i :\bbK_i \rightarrow \bbK_i$ be a continuous map
which is analytic on the interior, and
let $\A_i$ be a conjugacy algebra for $(\bbK_i, \eta_i)$ for $i = 1,2$.
Assume that the Fourier series map $\Delta_i$ is injective on $\A_i(A(\bbK_i))$.
Then any isomorphism $\gamma : \A_1(A(\bbK_1)) \rightarrow \A_2(A(\bbK_2)$
is automatically continuous.
\end{cor}

\Prf If $\S$ is any subset of $\A_2(A(\bbK_2)$,
\[ \bigcap_{n\ge0} \S U^n \subset \ker \Delta = \{0\} . \]
Thus if $\S(\gamma) \neq \{0\}$, then taking
$B_i = U$ in Lemma~\ref{Sinclair}, we obtain an integer $n_0$ so that
\[ \S(\gamma) U^{n_0} = \bigcap_{n\ge0} \S U^n = \{0\} .\]
Since right multiplication by $U$ is injective, $\S(\gamma) = \{0\}$.
Therefore $\gamma$ is continuous.
\bx

\begin{rem}
This same argument shows that an isomorphism between semicrossed products
$C_0(\X_i) \times_{\eta_i} \bbZ^+$ is automatically continuous.
However, when applied to a more general topological conjugacy algebra, one gets somewhat less.
One can conclude that $\S(\phi)$ is contained in $\bigcap_{n\ge1} \A U^n = \ker \Delta$.
In order to follow the plan of attack used here, we would need to characterize
$\ker \Delta$ algebraically.  When the set of fixed points has no interior, the
intersection $\J_n$ of the maximal ideals $\{ \theta_{x,0} : x \ne \eta(x) \}$ equals $\ker E_0$.
This set is determined by the fact that these functionals are not contained in the
closure of any analytic disk of characters. When the fixed point set has interior,
one also can try to use the ideal $\J_f$ obtained as the intersection of the kernels
of all maximal analytic disks of characters.  Unfortunately we do not see how to
recover $\ker \Delta$ from these two ideals.
\end{rem}

The following lemma establishes an important connection between nest
representations and differentiation.

\begin{lem} \label{derivative}
Let $\eta:\bbD \rightarrow \bbD$ be a continuous map which is analytic on the interior, and let
$\A$ be a topological conjugacy algebra for $(\bbD, \eta)$. Assume that there exists
$x \in \bbD^{\circ}$ so that $\eta(x) = x$. If $\pi \in \rep_{x, x} \A(A(\bbD))$
is a
continuous representation, then
\[
 \theta_{\pi, 1}(U)= \eta^{\prime}(x)\theta_{\pi, 2}(U).
\]
\end{lem}

\Prf Clearly, $\pi(I)$ is the identity matrix.
Let $\bz$ denote the identity function $\bz(z)= z$ for $z \in \bbD$.
Then there are scalars $a, w_i \in \bbC$ so that
\[
 \pi(\bz) =
 \left(\begin{array}{cc} x & a \\ 0 & x \end{array} \right) \qand
 \pi(U) =
  \left(\begin{array}{cc} w_1 & w_2 \\ 0 & w_3 \end{array} \right).
\]
Since the range of $\pi$ is not commutative, $a \neq 0$.  Therefore,
\[
 \pi(\bz^n) =
 \left(\begin{array}{cc} x^n & anx^{n-1} \\ 0 & x^n \end{array} \right)
\]
and by continuity,
\[
  \pi(f) =
  \left(\begin{array}{cc} f(x) & a f^{\prime}(x)
  \\ 0 & f(x) \end{array} \right) ,
\]
for any $f \in A(\bbD)$.

If we apply $\pi$ to the equation
$U f =  (f\circ \eta)U $, we obtain
\begin{equation*}
\left(\begin{array}{cc} w_1 & w_2 \\ 0 & w_3 \end{array} \right)
\left(\begin{array}{cc} f(x) & af'(x) \\ 0 & f(x)\end{array} \right)
=
\left(\begin{array}{cc} f(x) & af'(x)\eta'(x)\\ 0 & f(x) \end{array} \right)
\left(\begin{array}{cc} w_1 & w_2 \\ 0 & w_3 \end{array} \right)
\end{equation*}
By comparing the $(1,2)$-entries, we obtain
\[
aw_1f'(x) = a f'(x)\eta'(x)w_3 \qforal f \in A(\bbD),
\]
as desired.
\bx

\begin{lem} \label{rotation}
Consider the map $\eta_i(z)= c_iz$ for $z \in \bbD$,
where $|c_i|=1$ and let $\A_i$ be a conjugacy algebra for $(\bbD,\eta_i)$, $i=1,2$.
If the algebras $\A_1(A(\bbD))$ and $\A_2(A(\bbD))$ are isomorphic,
then $c_1=c_2$ or $c_1=c_2^{-1}$.
Furthermore, both possibilities occur for
the semicrossed products $A(\bbD) \times_{\eta_1}\bbZ^{+}$.
\end{lem}

\Prf  Assume that $\gamma: \A_1(A(\bbD))\to \A_2(A(\bbD))$ is an algebraic isomorphism.
By factoring out $\ker \Delta$ as in Corollary~\ref{kerDeltacor},
we may assume that $\Delta$ is injective on both $\A_i$.
Then by Corollary~\ref{auto_cont}, we see that $\gamma$ is continuous.

Let $\gamma_c: \M_{\A_1(A(\bbD))} \rightarrow \M_{\A_2(A(\bbD))}$ and
$\gamma_r: \rep_{\fT_2}  \A_1(A(\bbD)) \rightarrow$ \break $ \rep_{\fT_2} \A_2(A(\bbD))$
be defined as in the proof of Theorem \ref{main}.
Then $\gamma_r$ maps continuous representations
to continuous representations.

Once again, we are faced with two possibilities: either
$\gamma_c( \M_{\A_1(A(\bbD)), 0}) = \M_{\A_2(A(\bbD)), 0}$ or $\gamma_c( \M_{\A_1(A(\bbD)), 0}) =  \M_{\A_2(A(\bbD))}^0$.
In the first case we can proceed as in the proof of Theorem \ref{secmain} to conclude that
$c_1=c_2$. We therefore examine only the second case.

Since $\gamma_c( \M_{\A_1(A(\bbD)), 0}) =  \M_{\A_2(A(\bbD))}^0$, for each $z \in \bbD$ there exists
$\gamma_s(z) \in \bbD$ so that
\begin{equation} \label{strange}
\gamma_c(\theta_{z,0}) =\theta_{0,\gamma_s(x)}.
\end{equation}
Clearly, $\gamma_s$ is a bijection.
 From (\ref{strange}), it follows that
$\gamma_s(z) =\theta_{z,0}(\gamma^{-1}(U))$,
which is easily seen to imply that
$\gamma_s: \bbD\rightarrow \bbD$ is a homeomorphism
which is analytic on the interior.

Consider $z \in \bbD$, $z \neq0$, and
let $\rho \in \rep_{z, \eta_1(z)} \A_1(A(\bbD))$ be as in Example~\ref{example1}. Clearly,
\begin{equation} \label{rho}
\theta_{\gamma_r(\rho), 1}=\theta_{0,\gamma_s(z)} \text{ and }
\theta_{\gamma_r(\rho), 2}=\theta_{0,\gamma_s(\eta_1(z))}
\end{equation}
On the other hand, $\rho$ is continuous and therefore $\gamma_r(\rho)$
is continuous as well. Therefore, by applying Lemma \ref{derivative}
to (\ref{rho}), we obtain $\gamma_s(z) =\eta_2^{\prime}(0)\gamma_s(\eta_1(z))$,
i.e., $\gamma_s(z) = c_2\gamma_s (c_1z)$, for all $z \in \bbD\backslash \{0\}$.
This implies that $c_1=c_2^{-1}$.

For the last assertion, we need only examine the case where $c_2=c_1^{-1}$.
By \cite[Theorem II.4]{BP}, the semicrossed product
$A(\bbD) \times_{\eta_1}\bbZ^{+}$ is isomorphic
to the norm closed operator algebra generated
by a universal pair $(V, W)$ of unitary operators
satisfying $VW =c_1 WV$.
But then the pair $\{W,V\}$ is a universal pair of unitary operators
such that $WV = c_1^{-1}VW$.
So $A(\bbD) \times_{\eta_1}\bbZ^{+}$ is isomorphic
to $A(\bbD) \times_{\eta_1^{-1}}\bbZ^{+}$ as claimed.
\bx

\begin{thm} \label{maindisc2}
Let $\eta_i: \bbD \rightarrow \bbD$ be a continuous map which is analytic on the interior, $i=1,2$.
Assume further that $\eta_1$ is elliptic. Then:
\begin{enumerate}
\item If $\eta_1$ is a homeomorphism, then there are
isomorphic conjugacy algebras $\A_i(A(\bbK_i))$ for $(\bbK_i, \eta_i)$,
$i = 1,2$, if and only if
there exists a homeomorphism $\gamma_s: \bbD \rightarrow \bbD$
which is analytic on the interior and satisfies
$\gamma_s\circ \eta_1 = \eta_2\circ\gamma_s$ or
$\gamma_s\circ \eta_1 = \eta_2^{-1}\circ\gamma_s$.

\item If $\eta_1$ is not a homeomorphism, then there are
isomorphic conjugacy algebras $\A_i(A(\bbK_i))$ for $(\bbK_i, \eta_i)$,
$i = 1,2$, if and only if
there exists a homeomorphism $\gamma_s: \bbD \rightarrow \bbD$
which is analytic on the interior and satisfies
$\gamma_s\circ \eta_1 = \eta_2\circ\gamma_s$.
\end{enumerate}
\end{thm}

\Prf As in Lemma~\ref{rotation}, we may assume that $\Delta$ is injective
on both $\A_i$ and that the isomorphism $\gamma: \A_1(A(\bbD))\to \A_2(A(\bbD))$ is continuous.
Let $\gamma_c: \M_{\A_1(A(\bbD))} \to \M_{\A_2(A(\bbD))}$ and
$\gamma_r: \rep_{\fT_2}  \A_1(A(\bbD)) \to \rep_{\fT_2} \A_2(A(\bbD))$
be as in (\ref{induced1}) and (\ref{induced2}) respectively.
Note that by Corollary~\ref{maindisc1}, both $\eta_1$ and $\eta_2$ are elliptic.
Without loss of generality, we may assume that $\eta_1(0) = \eta_2(0)=0$.

We will show that if $\gamma_c( \M_{\A_1(A(\bbD)), 0}) =  \M_{\A_2(A(\bbD))}^0$, then
both $\eta_1$ and $\eta_2$ are homeomorphisms.
Proceeding as in the proof of Lemma \ref{rotation},
we obtain a homeomorphism $\gamma_s:\bbD \to \bbD$,
which is analytic on the interior and satisfies
$\gamma_s(z) =\eta_2^{\prime}(0)\gamma_s(\eta_1(z))$, $z \in \bbD\backslash\{0\}$.
 From this it follows that
$|\eta_2^{\prime}(0)|=1$, for otherwise $\gamma_s$ is not surjective.
By reversing our argument we conclude that
$|\eta_1^{\prime}(0)|=1$ as well.
The conclusion now follows from Schwarz's Lemma.

So if $\eta_1$ is not a homeomorphism, then
$\gamma_c( \M_{\A_1(A(\bbD)), 0}) = \M_{\A_2(A(\bbD)), 0}$.
So we may argue as in the proof of Theorem \ref{secmain}
to obtain the desired map $\gamma$. Conversely, if there exists
a conformal mapping $\gamma_s: \bbD \to \bbD$ so that
$\gamma_s\circ \eta_1 = \eta_2\circ\gamma_s$,
then the two semicrossed products are easily seen to be isomorphic.
This proves (ii).

Finally assume that $\eta_1$ is a homeomorphism.
Then, by (ii) $\eta_2$ is also a homeomorphism.
Any elliptic conformal mapping of the unit disc is conformally
conjugate to a map fixing the origin, namely a rotation.
Thus the conclusion follows from Lemma~\ref{rotation}.
\bx

Applying the main result to the case of the
analytic semi-crossed products, we obtain:

\begin{cor}\label{analyticsemicross}
Let $\bbK_i \in \K$, and
let $\eta_i :\bbK_i \rightarrow \bbK_i$ be continuous maps
which are analytic on the interior for $i = 1,2$.
The algebras $A(\bbK_1) \times_{\eta_1} \bbZ^+$
and $A(\bbK_2) \times_{\eta_1} \bbZ^+$
are isomorphic if and only if either
\begin{enumerate}
\item $\eta_1$ and $\eta_2$ are analytically conjugate; or

\item $\bbK_i$ are simply connected,
$\eta_i$ are homeomorphisms with a unique fixed point,
this point is in the interior of $\bbK_i$,
and $\eta_2$ is analytically conjugate to $\eta_1^{-1}$.
\end{enumerate}
\end{cor}

\begin{rem}
It is certainly possible for analytic maps of the unit disk
to be topologically conjugate but not analytically conjugate.
A simple example is given by $\eta_1(z) = z/2$ and $\eta_2(z) = z/4$.
Then $\gamma(re^{i\theta}) = r^2e^{i\theta}$ satisfies
$\gamma \circ \eta_1 = \eta_2 \circ \gamma$.
Schwarz's Lemma shows that no conformal automorphism of the disk
can carry the disk of radius $1/2$ onto the disk of radius $1/4$.
So they are not analytically conjugate.

This can even occur for conformal automorphisms of the disk.
Let $\eta_1(z) = \dfrac{z-1/2}{1-z/2}$ and
$\eta_2(z) = \dfrac{z-1/4}{1-z/4}$.
If one identifies the disk with the left half plane via the
conformal map $w=\dfrac{z-1}{z+1}$ that sends $1$ to $0$
and $-1$ to the point at infinity,
then this conjugates both maps to dilations
$\tau_1(w) = w/3$ and $\tau_2(w) = 3w/5$.
Again these maps are topologically conjugate.

Recall that the conformal automorphisms of
the disk are just the Mobius maps.  Since $\eta_i$ both
fix the points $\pm1$, with $1$ being the unique contractive
fixed point, the intertwining homeomorphism $\gamma$ must
fix $\pm1$.  But all of these maps commute with $\eta_i$.
So they are not analytically conjugate.
\end{rem}

%%%%%%%%%%%%%%%%%%%%%%%%%%%%%%%%%%%%%%%%%%%%%

\end{document}